\documentclass[fleqn,a4paper,12pt]{article}
\usepackage{amsmath,graphicx,amssymb,amsthm,amscd}

\newcommand{\diag}{{\rm diag}}

\newcommand{\bfA}{\mbox{\boldmath $A$}}

\newcommand{\bfE}{\mbox{\boldmath $E$}}
\newcommand{\bfG}{\mbox{\boldmath $G$}}

\newcommand{\bfI}{\mbox{\boldmath $I$}}
\newcommand{\bfJ}{\mbox{\boldmath $J$}}

\newcommand{\bfO}{\mbox{\boldmath $O$}}

\newcommand{\bfQ}{\mbox{\boldmath $Q$}}

\newcommand{\bfdelta}{\mbox{\boldmath $\delta$}}

\newcommand{\bfz}{\mbox{\boldmath $0$}}
\newcommand{\bfi}{\mbox{\boldmath $1$}}

\newcommand{\bfq}{\mbox{\boldmath $q$}}
\newcommand{\bfu}{\mbox{\boldmath $u$}}

\newcommand{\bfv}{\mbox{\boldmath $v$}}

\newtheorem{lemma}{Lemma}
\newtheorem{prop}{Proposition}
\newtheorem{theorem}{Theorem}
\newtheorem{zu}{Fig.}

\makeatletter
\@addtoreset{equation}{section}%
\@addtoreset{lemma}{section}%
\@addtoreset{prop}{section}%
\@addtoreset{theorem}{section}%
\@addtoreset{zu}{section}%

\makeatother

\begin{document}
\title{The best constant of discrete Sobolev inequality on 
the C60 fullerene buckyball}

\author{
Yoshinori Kametaka$^1$, 
Atsushi Nagai$^2$, 
Hiroyuki Yamagishi$^3$, \\
Kazuo Takemura$^4$ and 
Kohtaro Watanabe$^5$
}
\maketitle


\begin{abstract}
The best constants of two kinds of discrete Sobolev inequalities on 
the C60 fullerene buckyball are obtained. 
All the eigenvalues of discrete Laplacian $\bfA$ corresponding to
the buckyball are found. They are roots of 
algebraic equation at most degree $4$ with integer coefficients.
Green matrix $\bfG(a)=(\bfA+a\bfI)^{-1}\ (0<a<\infty)$ and 
the pseudo Green matrix $\bfG_*=\bfA^{\dagger}$ 
are obtained by using computer software Mathematica. 
Diagonal values of $\bfG_*$ and $\bfG(a)$ are identical and
they are equal to the best constants of discrete Sobolev inequalities. 

\bigskip



\end{abstract}
\section{Introduction}
The best constant of discrete Sobolev inequality was obtained for 
Platon's regular polyhedra \cite{rph,rph2} 
and truncated regular tetra-, hexa- and octa- polyhedra \cite{trph}. 
The purpose of this paper is to extend the above result to 
the truncated regular icosahedron, that is, the C60 Fullerene Buckyball \cite{c60}. 
The studies of discrete Sobolev inequalities start with \cite{dsob,dber},
where inequalities on a periodic lattice are investigated. 

\begin{figure}[htbp]
\includegraphics[scale=0.5]{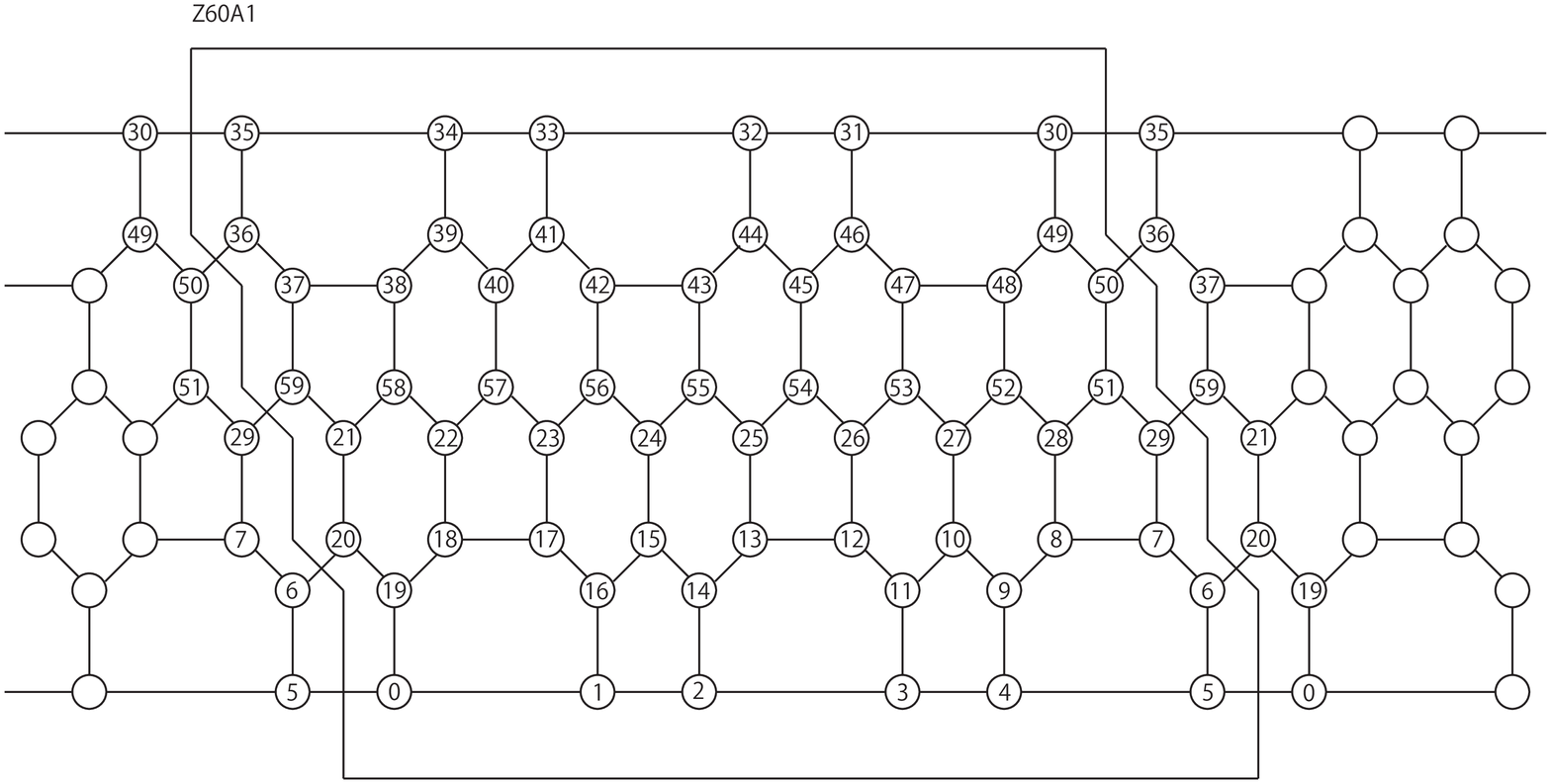}
\begin{zu}\label{a1z}The zigzag type of C60 fullerene Buckyball.\end{zu}
\end{figure}

\begin{figure}[htbp]
\includegraphics[scale=0.5]{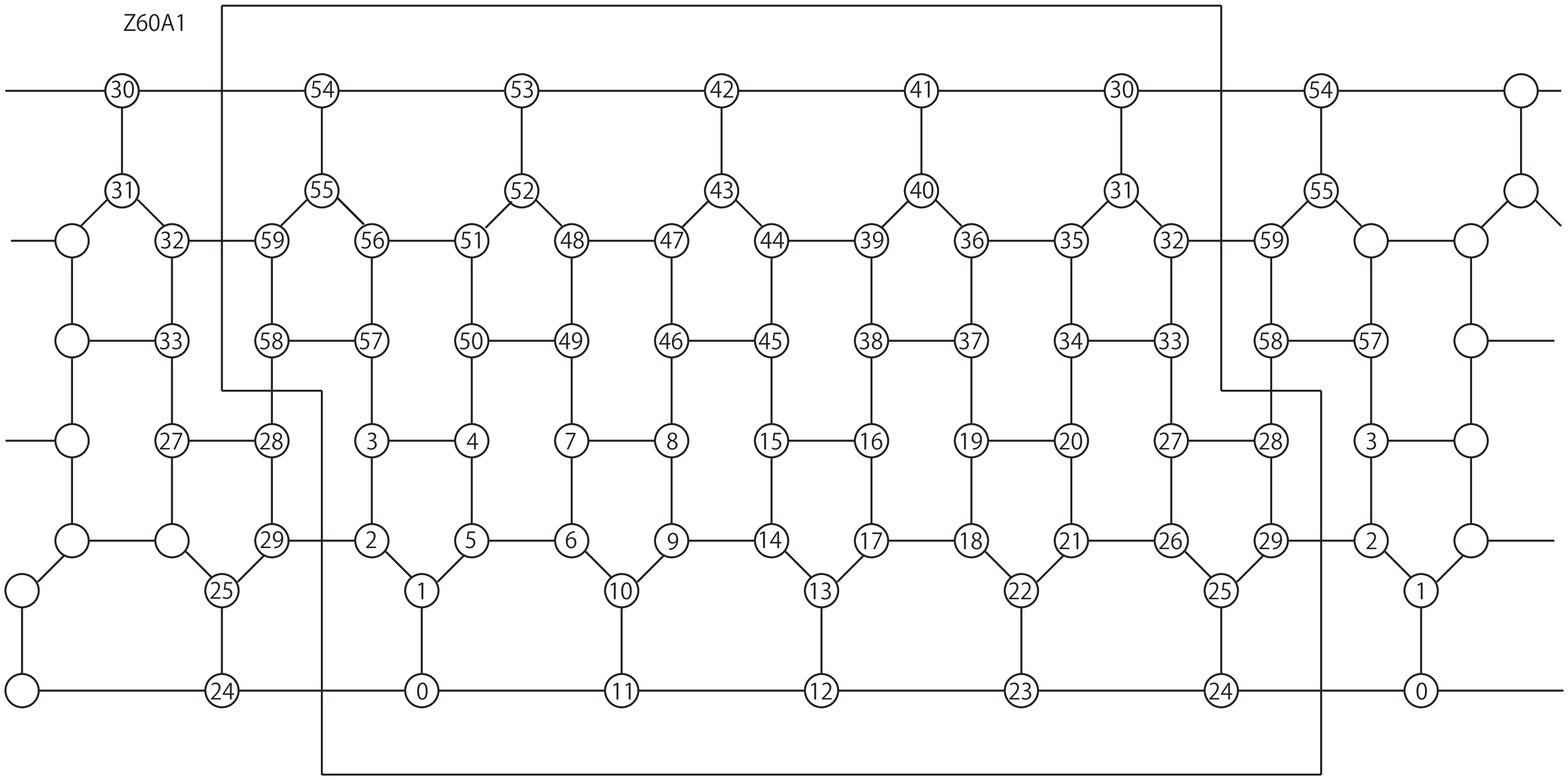}
\begin{zu}\label{a1a}The armchair type of C60 fullerene Buckyball.\end{zu}
\end{figure}

Buckyball, or BB for short, is a polyhedron which consists of 
60 carbon atoms, 
including 12 pentagons and 20 hexagons. 
The numbers of vertices $v$, faces $f$ and edges $e$ are $v=60, f=32, e=90$, 
which satisfies Euler polyhedron theorem $v+f=e+2$. 
We use the following two equivalent simplest plane graph expressions 
in Figs.\ref{a1z} and \ref{a1a}. 
Fig.\ref{a1z} is called the zigzag type and 
Fig.\ref{a1a} is called the armchair type of BB.

This paper is organized as follows; 
In section 2, we introduce discrete Laplacian corresponding to BB. 
In section 3, we present main theorems in this paper. 
In section 4, we investigate pseudo Green matrix $\bfG_*$ and Green matrix $\bfG(a)$. 
In section 5, we give the spectral decomposition of $\bfA$ and investigate diagonal
values of matrices $\bfG_*$ and $\bfG(a)$.
In section 6 and 7, we prove main theorems. 

\section{Discrete Laplacian}
\begin{figure}[hbtp]
\begin{minipage}[cbt]{6cm}
{\includegraphics[scale=1.5]{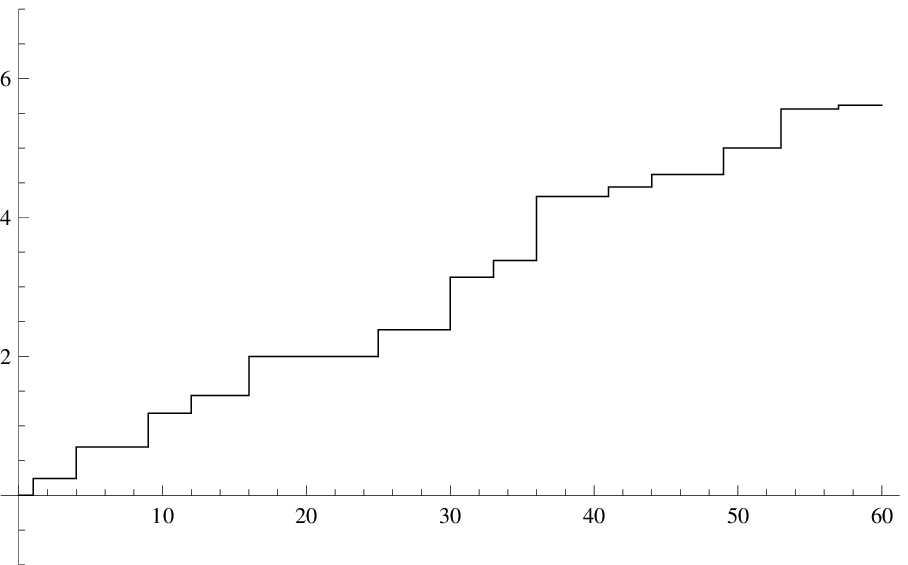}}
\end{minipage}
\begin{zu}\label{eigentd}Eigenvalue distribution of C60 fullerene buckyball.\end{zu}
\end{figure}

We treat a classical mechanical model of BB. 
Each neighboring two atoms are connected by a classical linear spring 
with uniform spring constant. 
The vertices are labeled by integers $i\ (0\le i\le 59)$ as in Figs. \ref{a1z} and \ref{a1a}. 
The edge set $e$ is a set of $(i,j)$, where $i$ and $j$ are connected by an edge. 
Discrete Laplacian $\bfA$ is defined by
\begin{align*}
&\bfA=\Big(a(i,j)\Big)\quad (0\le i,j\le 59)
,\qquad
a(i,j)=
\left\{\begin{array}{ll}
3\ &(i=j),\\  -1\ &((i,j)\in e),\\ 0\ &(\mbox{else}).
\end{array}\right.
\end{align*}
$\bfA$ is a $60\times 60$ real symmetric 
positive-semidefinite matrix. 
Although two discrete Laplacians derived from Figs.~\ref{a1z} and \ref{a1a} are different, they are similar,
they possess a common characteristic polynomial $P(x)$ given by 
\begin{align*}
&P(x)=\det(x\bfI-\bfA)=\\&
x(x-2)^9(x-5)^4(x^2-5x+3)^5(x^2-7x+11)^5\times\\
&(x^2-7x+8)^4(x^2-9x+19)^3
(x^4-9x^3+25x^2-22x+4)^3
.
\end{align*}
Each factors of $P(x)$ are a polynomial of $x$ with integer coefficients 
with degree at most $4$. 
Eigenvalues of $\bfA$, their multiplicities and approximate values, are shown in 
Fig.~\ref{eigent}.
Fig.~\ref{eigentd} illustrates a distribution of 
eigenvalues $\lambda_j$ of $\bfA$. 
The eigenvector corresponding to $\lambda_0=0$ is $\bfi={}^t(1,\cdots,1)\in {\bf C}^{60}$. 

\begin{figure}[hbtp]
\begin{minipage}[cbt]{5cm}
\begin{tabular}{|c|c|c|c|}
\hline
number &eigenvalue&approximate&multiplicity \\ 
$j$&$\lambda_j$&value&\\
\hline 
0&0&0&1\\
 1$\sim$ 3&$(9 - \sqrt{5} - \sqrt{38 - 2 \sqrt{5}})/4$&0.24&3\\
 4$\sim$ 8&$(5 - \sqrt{13})/2$&0.69&5\\
 9$\sim$11&$(9 + \sqrt{5} - \sqrt{38 + 2 \sqrt{5}})/4$&1.17&3\\
12$\sim$15&$(7 - \sqrt{17})/2$&1.43&4\\
16$\sim$24&$2$&2&9\\
25$\sim$29&$(7 - \sqrt{5})/2$&2.38&5\\
30$\sim$32&$(9 - \sqrt{5} + \sqrt{38 - 2 \sqrt{5}})/4$&3.13&3\\
33$\sim$35&$(9 - \sqrt{5})/2$&3.38&3\\
36$\sim$40&$(5 + \sqrt{13})/2$&4.30&5\\
41$\sim$43&$(9 + \sqrt{5} + \sqrt{38 + 2 \sqrt{5}})/4$&4.43&3\\
44$\sim$48&$(7 + \sqrt{5})/2$&4.61&5\\
49$\sim$52&$5$&5&4\\
53$\sim$56&$(7 + \sqrt{17})/2$&5.56&4\\
57$\sim$59&$(9 + \sqrt{5})/2$&5.61&3\\
\hline
\end{tabular}
\end{minipage}
\begin{zu}\label{eigent}Eigenvalues of C60 fullerene buckyball.\end{zu}
\end{figure}

\section{Conclusion}
For each vertex $i$, we attach a complex number $u(i)$. 
For vector $\bfu={}^t(u(0),\cdots,u(59))\in {\bf C}^{60}$, 
we introduce two kinds of Sobolev energies, that is, potential energies. 
\begin{align*}
&E(\bfu)=\sum_{(i,j)\in e}|u(i)-u(j)|^2=\bfu^{*}\bfA\bfu,
\\
&E(a,\bfu)=E(\bfu)+a\sum_{j=0}^{59}\!|u(j)|^2=\bfu^{*}(\bfA+a\bfI)\bfu.
\end{align*}
The constant $a\ (0<a<\infty)$ stands for a dumping parameter. 
The relations 
\begin{align}\label{gs-def}
\bfA\bfG_*=\bfG_*\bfA=\bfI-\bfE_0,\qquad 
\bfG_*\bfE_0=\bfE_0\bfG_*=\bfI-\bfE_0
\end{align}
determine a unique solution $\bfG_*$. 
$\bfE_0=(1/60)\bfi{}^t\bfi$ is 
the orthogonal projection matrix to the eigenspace corresponding to eigenvalue 0 of $\bfA$.  
We call $\bfG_*$ the pseudo Green matrix of $\bfA$. 
$\bfG_*$ coincides with a Penrose-Moore generalized inverse matrix $\bfA^{\dagger}$ of $\bfA$. 
$\bfG(a)=(\bfA+a\bfI)^{-1}\ (0<a<\infty)$ is Green matrix of $\bfA$. 
In section 4, we show 
\begin{align}\label{gs}
\bfG_*=\lim_{a\to +0}\left(\bfG(a)-a^{-1}\bfE_0\right)
. 
\end{align}

The important fact is that diagonal elements of $\bfG_*$ 
and $\bfG(a)$ are all the same, which is observed by direct calculations of 
Mathematica. 
Concrete values of these are shown by \eqref{C0} and \eqref{C0a} into 
following two theorems. 
The above fact reflects high symmetry of BB. 
\begin{theorem}\label{th1}
For any $\bfu\in {\bf C}^{60}$ satisfying $u(0)+\cdots+u(59)=0$, 
there exists a positive constant $C$ which is independent of $\bfu$, 
such that the discrete Sobolev inequality
\begin{align}\label{dsob1}
\left(\max_{0\le j\le 59}|\,u(j)\,|\right)^2\le C E(\bfu)
\end{align}
holds. Among such $C$, the best (the smallest) constant $C_0$ is
\begin{align}\label{C0}
&C_0=
\frac{239741}{376200}
\fallingdotseq
0.63727\cdots
.
\end{align}
Diagonal values of $\bfG_*$ take the same value $C_0$. 
In particular, we have 
\begin{align}\label{C0lam}
C_0=\frac{1}{60}\sum_{j=1}^{59}\frac{1}{\lambda_j}.
\end{align}
If we replace $C$ by $C_0$ in the above inequality, 
then the equality holds for any column vector $\bfu$ of $\bfG_*$. 
\end{theorem}

\begin{theorem}\label{th2}
For any $\bfu\in {\bf C}^{60}$, 
there exists a positive constant $C$, which is independent of $\bfu$, 
such that the discrete Sobolev inequality
\begin{align}\label{dsob2}
\left(\max_{0\le j\le 59}|\,u(j)\,|\right)^2\le C E(a,\bfu)
\end{align}
holds. Among such $C$, the best constant $C(a)$ is 
\begin{align}\label{C0a}
&C(a)=\frac{N(a)}{D(a)},
\\&
N(a)=
3344+160806a+1153562a^2+3594661a^3+
6334271a^4+7104785a^5+
\nonumber\\&
5406109a^6+2893077a^7+
1109403a^{8}+306415a^{9}+60463a^{10}+8315a^{11}+
\nonumber\\&757a^{12}+41a^{13}+a^{14}
,
\nonumber\\&
D(a)=
a(2+a)(5+a)(3+5a+a^2)(8+7a+a^2)\times
\nonumber\\&
(11+7a+a^2)(19+9a+a^2)(4+22a+25a^2+9a^3+a^4)
\nonumber
.
\end{align}
Diagonal values of $\bfG(a)$ take the same value $C(a)$. 
In particular, we have 
\begin{align}\label{C0alam}
&C(a)=\frac{1}{60}\sum_{j=0}^{59}\frac{1}{\lambda_j+a}.
\end{align}
If we replace $C$ by $C(a)$ in the above inequality, 
the equality holds for any column vector $\bfu$ of $\bfG(a)$. 
\end{theorem}

\begin{figure}[hbtp]
\begin{minipage}[cbt]{7.5cm}
\includegraphics[scale=0.8]{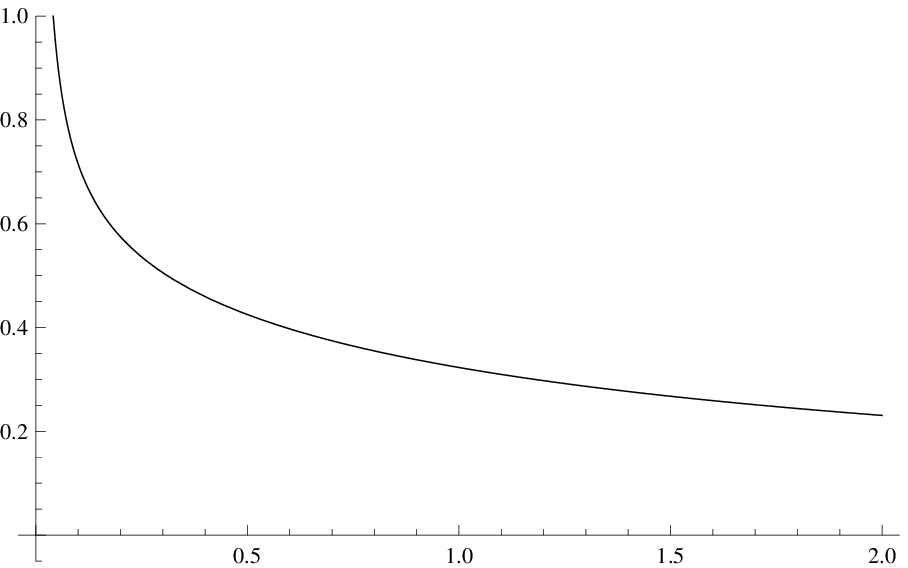}
\vspace{-3mm}\begin{zu}\label{z60gha}The graph of $C(a)$.\end{zu}
\end{minipage}
%
\begin{minipage}[cbt]{7.5cm}
\includegraphics[scale=0.8]{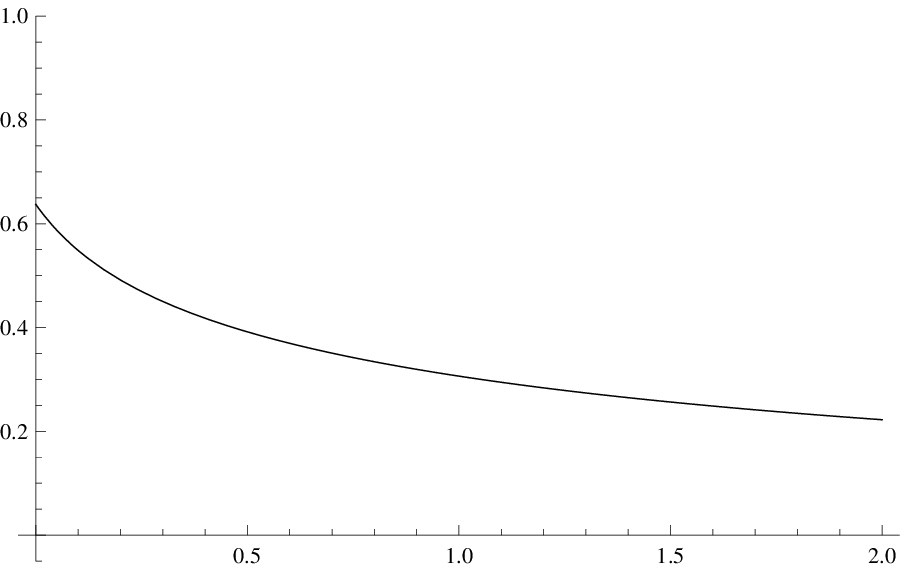}
\vspace{-3mm}\begin{zu}\label{z60ghgs}The graph of $C(a)-(60a)^{-1}$.\end{zu}
\end{minipage}
\end{figure}

From \eqref{C0alam}, it follows that $C(a)$ is a monotone decreasing function of $a$. 
Fig.\ref{z60gha} is the graph of $C(a)$ and 
Fig.\ref{z60ghgs} is the graph of $C(a)-(60a)^{-1}$. 
From \eqref{gs}, \eqref{C0lam} and \eqref{C0alam}, we have the following theorem. 
\begin{theorem}\label{th3}
The relation 
\begin{align}\label{C0lim}
C_0=\lim_{a\to 0}\left(C(a)-(60a)^{-1}\right)
\end{align}
holds. 
\end{theorem}

Physical meaning of the above theorems is as follows.
If $u(i)$ takes a real value, then $u(i)$ represents 
a deviation from the steady state. The discrete Sobolev inequalities shows that
the square of maximum of deviation $|u(i)|$ is estimated from above by 
constant multiples of the potential energies $E(\bfu)$ and $E(a,\bfu)$.  
Hence, it is expected that the best constants $C_0$ and $C(a)$ represent rigidity 
of the mechanical model.

\section{Green matrix and pseudo Green matrix}
The discrete Laplacian $\bfA$ is written in the following form
\begin{align*}
&\bfA=
\left(\begin{array}{c|c}
\bfA_0 & \bfA_1 \\ \hline
\bfA_1 & \bfA_0
\end{array}\right),
\end{align*}
where $\bfA_0$ and $\bfA_1$ are $30\times 30$ real symmetric matrices. 
We introduce a $60\times 60$ matrix $\bfJ$ given by
\begin{align*}
\bfJ=\left(\begin{array}{c|c}
\bfI & \bfI \\ \hline
\bfI & -\bfI
\end{array}\right),
\end{align*}
where $\bfI$ is a $30\times 30$ identity matrix. We have $\bfJ^{-1}=(1/2)\bfJ$. 
Using $\bfJ$, we have
\begin{align*}
\bfJ^{-1}\bfA\bfJ=
\left(\begin{array}{c|c}
\bfA_{+} & \bfO \\ \hline
\bfO & \bfA_{-}
\end{array}\right)
,\qquad 
\bfA_{\pm}=\bfA_0\pm \bfA_1
.
\end{align*}
$\bfA_+$ has an eigenvalue 0, corresponding eigenvector 
is $\bfi={}^t(1,\cdots,1)\in {\bf C}^{30}$. 
The other eigenvalues of $\bfA_+$ and $\bfA_-$ are all positive. 
%
%
The pseudo Green matrix $\bfG_*$ of $\bfA$ is 
\begin{align}\label{eq:bfG}
&\bfG_*=
\left(\begin{array}{c|c}
\bfG_{*0} & \bfG_{*1} \\ \hline
\bfG_{*1} & \bfG_{*0}
\end{array}\right)
=\bfJ
\left(\begin{array}{c|c}
\bfG_{*+} & \bf0 \\ \hline
\bf0 & \bfG_{*-}
\end{array}\right)
\bfJ^{-1},
\end{align}
where $\bfG_{*+}=\bfA_{+}^{\dagger}$ is a Penrose-Moore generalized inverse matrix of $\bfA_{+}$ 
and $\bfG_{*-}=\bfA_{-}^{-1}$ is an inverse matrix of $\bfA_{-}$. 
We have
\begin{align*}
&\bfG_{*\pm}=\bfG_{*0}\pm \bfG_{*1}
\Leftrightarrow
\left\{\begin{array}{ll}
\bfG_{*0}=(\bfG_{*+}+\bfG_{*-})/2,
\\
\bfG_{*1}=(\bfG_{*+}-\bfG_{*-})/2.
\end{array}\right.
\end{align*}
It is easy to see that
\begin{align*}
\bfA+a\bfI=
\bfJ\left(\begin{array}{c|c}
\bfA_{+}+a\bfI & \bfO \\ \hline
\bfO & \bfA_{-}+a\bfI
\end{array}\right)
\bfJ^{-1}. 
\end{align*}
If we put $\bfG_{\pm}(a)=(\bfA_{\pm}+a\bfI)^{-1}$, 
Green matrix $\bfG(a)=(\bfA+a\bfI)^{-1}$ is given by
\begin{align}\label{eq:bfGa}
\bfG(a)=(\bfA+a\bfI)^{-1}=
\bfJ\left(\begin{array}{c|c}
\bfG_{+} & \bfO \\ \hline
\bfO & \bfG_{-}
\end{array}\right)\!\!(a)
\bfJ^{-1}
=
\left(\begin{array}{c|c}
\bfG_{0} & \bfG_{1} \\ \hline
\bfG_{1} & \bfG_{0}
\end{array}\right)\!\!(a),
\end{align}
where
\begin{align*}
&\bfG_{\pm}(a)=\bfG_0(a)\pm \bfG_1(a)
\Leftrightarrow
\left\{\begin{array}{ll}
\bfG_0(a)=(\bfG_+(a)+\bfG_-(a))/2,
\\
\bfG_1(a)=(\bfG_+(a)-\bfG_-(a))/2.
\end{array}\right.
\end{align*}
Using the above expressions \eqref{eq:bfG} anf \eqref{eq:bfGa}, one can reduce
the number of computations in finding explicit forms of $\bfG_*$ and $\bfG(a)$.
\section{The spectral decomposition of $\bfA$}
%
For $0\le j\le 59$, we introduce 
\begin{align*}
\bfdelta_{j}={}^t(\delta(-j),\delta(1-j),\cdots,\delta(59-j))
,\qquad 
\delta(k)=\left\{\begin{array}{ll}
1 ~ &({\rm Mod}(k,60)=0),\\ 0 ~&({\rm Mod}(k,60)\ne 0).
\end{array}\right.
\end{align*}
The discrete Laplacian $\bfA$ is a $60\times 60$ real symmetric matrix. 
Using unitary matrix $\bfQ$, $\bfA$ is diagonalized as 
$\bfQ^{*}\bfA\bfQ=\widetilde{\bfA}=\diag\{\lambda_0,\cdots,\lambda_{59}\}$. 
$\bfE_k=\bfq_k\bfq_k^*\ (0\le k\le 59)$ are orthogonal projection matrices 
which satisfy
\begin{align*}
\bfE_k^*=\bfE_k,\quad
\bfE_k\bfE_l=\delta(k-l)\bfE_k\qquad (0\le k,l\le 59). 
\end{align*}
The spectral decomposition of 
$60\times 60$ unit matrix $\bfI$ and the discrete Laplacian $\bfA$ are rewritten as
\begin{align*}
&\bfI=\bfQ\bfQ^*=\sum_{k=0}^{59}\bfE_k,
\qquad 
\bfA=\bfQ\widetilde{\bfA}\bfQ^*=
\sum_{k=0}^{59}\lambda_k\bfE_k=
\sum_{k=1}^{59}\lambda_k\bfE_k
.
\end{align*}
For $0<a<\infty$, 
we have
\begin{align}
\label{sp-g}
&\bfG(a)=(\bfA+a\bfI)^{-1}=
\sum_{k=0}^{59}(\lambda_k+a)^{-1}\bfE_k
\\
\label{sp-gs}
&\bfG_*=\lim_{a\to +0}\left(\bfG(a)-a^{-1}\bfE_0\right)=\sum_{k=1}^{59}\lambda_k^{-1}\bfE_k.
\end{align}
It is easy to see that \eqref{sp-gs} satisfies \eqref{gs-def}. 

The rest of this section is devoted to the proof of \eqref{C0lam}
and \eqref{C0alam} in the main theorems.

\begin{prop}\label{prop-c0}
For any $j_0\ (0\le j_0\le 59)$, we have
\begin{align*}
&C_0={}^t\bfdelta_{j_0}\bfG_*\bfdelta_{j_0}=
\frac{1}{60}\sum_{k=1}^{59}\lambda_k^{-1}.
\end{align*}
\end{prop}

\noindent
{\bf Proof of Proposition \ref{prop-c0}}\quad
First of all, we have
\begin{align}\label{e}
&\sum_{j=0}^{59}{}^t\bfdelta_j\bfE_k\bfdelta_j=
\sum_{j=0}^{59}{}^t\bfdelta_j\bfq_k\bfq_k^{*}\bfdelta_j=
\sum_{j=0}^{59}|\bfq_k^*\bfdelta_j|^2=
\bfq_k^*\bfq_k=1
\quad (0\le k \le 59).
\end{align}
We note that diagonal values of $\bfG_*$ are identical. 
For any $j_0\ (0\le j_0\le 59)$, using \eqref{sp-gs} and \eqref{e}, we have
\begin{align*}
&C_0={}^t\bfdelta_{j_0}\bfG_*\bfdelta_{j_0}=
\frac{1}{60}\sum_{j=0}^{59}{}^t\bfdelta_j\bfG_*\bfdelta_j=
\frac{1}{60}\sum_{j=0}^{59}{}^t\bfdelta_j\sum_{k=1}^{59}\lambda_k^{-1}\bfE_k\bfdelta_j=
\\&
\frac{1}{60}\sum_{k=1}^{59}\lambda_k^{-1}
\sum_{j=0}^{59}{}^t\bfdelta_j\bfE_k\bfdelta_j=
\frac{1}{60}\sum_{k=1}^{59}\lambda_k^{-1}
.
\end{align*}
Thus we have Proposition \ref{prop-c0}. 
That is to say, we have \eqref{C0lam}. 
\hfill{$\blacksquare$}

\bigskip

\begin{prop}\label{prop-c0a}
For any $j_0\ (0\le j_0\le 59)$, we have
\begin{align*}
&C(a)={}^t\bfdelta_{j_0}\bfG(a)\bfdelta_{j_0}=
\frac{1}{60}\sum_{k=0}^{59}(\lambda_k+a)^{-1}.
\end{align*}
\end{prop}

\noindent
{\bf Proof of Proposition \ref{prop-c0a}}\quad
We note that the diagonal values of $\bfG(a)$ are identical. 
For any $j_0\ (0\le j_0\le 59)$, using \eqref{sp-g} and \eqref{e}, we have
\begin{align*}
&C(a)=
{}^t\bfdelta_{j_0}\bfG(a)\bfdelta_{j_0}=
\frac{1}{60}\sum_{j=0}^{59}{}^t\bfdelta_j\bfG(a)\bfdelta_j=
\frac{1}{60}\sum_{j=0}^{59}\bfdelta_j\sum_{k=0}^{59}{}^t(\lambda_k+a)^{-1}\bfE_k\bfdelta_j=
\\&
\frac{1}{60}\sum_{k=0}^{59}(\lambda_k+a)^{-1}
\sum_{j=0}^{59}{}^t\bfdelta_j\bfE_k\bfdelta_j=
\frac{1}{60}\sum_{k=0}^{59}(\lambda_k+a)^{-1}
.
\end{align*}
Thus we have Proposition \ref{prop-c0a}. 
That is to say, we have \eqref{C0alam}. 
\hfill{$\blacksquare$}

\section{Proof of Theorem \ref{th1}}
For $\bfu, \bfv\in {\bf C}^{60}$, we attach an inner product $(\bfu,\bfv)$ defined by 
\begin{align*}
(\bfu,\bfv)=\bfv^{*}\bfu,\qquad \|\,\bfu\,\|^2=(\bfu,\bfu).
\end{align*} 
For $\bfu, \bfv \in {\bf C}_0^{60}:= \{\bfu | \bfu\in {\bf C}^{60}\ and\ 
u(0)+\cdots+u(59)=0\}$, we also define $(\bfu,\bfv)_A$ by  
\begin{align*}
&(\bfu,\bfv)_A=(\bfA\bfu,\bfv)=\bfv^{*}\bfA\bfu,
\qquad 
\|\,\bfu\,\|_A^2=(\bfu,\bfu)_A=E(\bfu).
\end{align*}
\begin{lemma}\label{rep1}
For any $\bfu\in {\bf C}_0^{60}$, 
we have the following reproducing relation.
\begin{align}\label{rep10}
u(j)=
(\bfu,\bfG_*\bfdelta_j)_A
\qquad (0\le j\le 59). 
\end{align}
\end{lemma}

\noindent{\bf The proof of Lemma \ref{rep1}}\quad
We first note that $\bfG_*$ is an Hermitian matrix and that 
$\bfE_0\bfu=60^{-1}\bfi{}^t\bfi\bfu=\bfz$. 
For any $\bfu\in {\bf C}_0^{60}$ 
and any fixed $j=0,1,\cdots,59$, we have
\begin{align*}
&(\bfu,\bfG_*\bfdelta_j)_A=(\bfA\bfu,\bfG_*\bfdelta_j)=
{}^t\bfdelta_j\bfG_*\bfA\bfu=
{}^t\bfdelta_j(\bfI-\bfE_0)\bfu=
{}^t\bfdelta_j\bfu=
u(j)
.
\end{align*}
This completes the proof of Lemma \ref{rep1}. 
\hfill{$\blacksquare$}

\bigskip

Putting $\bfu=\bfG_*\bfdelta_j$ in \eqref{rep10}, we have
\begin{align}\label{rep11}
C_0={}^t\bfdelta_j\bfG_*\bfdelta_j
=\|\bfG_*\bfdelta_j\|_A^2=E(\bfG_*\bfdelta_j).
\end{align}
Applying Schwarz inequality to \eqref{rep10} and 
using \eqref{rep11}, we have
\begin{align*}
|\,u(j)\,|^2 \le \|\,\bfu\,\|_A^2 \|\,\bfG_*\bfdelta_j\,\|_A^2 =
C_0E(\bfu)
.
\end{align*}
Taking the maximum with respect to $j$ on both sides, we obtain 
discrete Sobolev inequality
\begin{align}\label{best-sob1}
\left(\max_{0\le j\le 59}|\,u(j)\,|\right)^2 \le C_0E(\bfu).
\end{align}
If we take $\bfu=\bfG_*\bfdelta_{j_0}$ in \eqref{best-sob1}, then we have
\begin{align*}
\left(\max_{0\le j\le 59}|\,{}^t\bfdelta_{j}\bfG_*\bfdelta_{j_0}\,|\right)^2\le
C_0 E(\bfG_*\bfdelta_{j_0}) = C_0^2
.
\end{align*}
Combining this and a trivial inequality
\begin{align*}
C_0^2=\left({}^t\bfdelta_{j_0}\bfG_*\bfdelta_{j_0}\right)^2
\le 
\left(\max_{0\le j\le 59}|\,{}^t\bfdelta_{j}\bfG_*\bfdelta_{j_0}\,|\right)^2,
\end{align*}
we have
\begin{align*}
\left(\max_{0\le j\le 59}|\,{}^t\bfdelta_{j}\bfG_*\bfdelta_{j_0}\,|\right)^2 = 
C_0 E(\bfG_*\bfdelta_{j_0})
.
\end{align*}
This shows that if we replace $C$ by $C_0$ in an inequality \eqref{dsob1}, 
the equality holds for $\bfG_*\bfdelta_{j_0}$. 
\eqref{C0} follows from \eqref{C0lam} by simple calculation. 
\eqref{C0} is also confirmed by the observation of the diagonal 
values of $\bfG_*$. 
This completes the proof of Theorem \ref{th1}.\hfill{$\blacksquare$}

\section{Proof of Theorem \ref{th2}}
For $\bfu, \bfv\in {\bf C}^{60}$, we introduce 
\begin{align*}
&(\bfu,\bfv)_H=((\bfA+a\bfI)\bfu,\bfv)=\bfv^{*}(\bfA+a\bfI)\bfu,
\qquad
\|\,\bfu\,\|_H^2=(\bfu,\bfu)_H=E(a,\bfu).
\end{align*}
In this section, we use the simpler notation $\bfG=\bfG(a)$. 

\begin{lemma}\label{rep2}
For any $\bfu\in {\bf C}^{60}$, we have the following reproducing relation.
\begin{align}\label{rep20}
u(j)=
(\bfu,\bfG\bfdelta_j)_H
\qquad (0\le j\le 59)
.
\end{align}
\end{lemma}

\noindent{\bf The proof of Lemma \ref{rep2}}\quad
We note that $\bfG$ is an Hermitian matrix. For any $\bfu\in {\bf C}^{60}$ 
and any fixed $j=0,1,\cdots,59$, we have
\begin{align*}
&(\bfu,\bfG\bfdelta_j)_H=
((\bfA+a\bfI)\bfu,\bfG\bfdelta_j)=
{}^t\bfdelta_j\bfG(\bfA+a\bfI)\bfu=
{}^t\bfdelta_j \bfu=u(j).
\end{align*}
This completes the proof of Lemma \ref{rep2}. 
\hfill{$\blacksquare$}

\bigskip

Putting $\bfu=\bfG\bfdelta_j$ in \eqref{rep20}, we have
\begin{align}\label{rep21}
C(a)={}^t\bfdelta_j\bfG\bfdelta_j=
\|\bfG\bfdelta_j\|_H^2=E(a,\bfG\bfdelta_j).
\end{align}
Applying Schwarz inequality to \eqref{rep20} and 
using \eqref{rep21}, we have
\begin{align*}
|\,u(j)\,|^2 \le \|\,\bfu\,\|_H^2 \|\,\bfG\bfdelta_j\,\|_H^2 = 
C(a)E(a,\bfu)
.
\end{align*}
Taking the maximum with respect to $j$ on both sides, we have the discrete Sobolev inequality
\begin{align}\label{best-sob2}
\left(\max_{0\le j\le 59}|\,u(j)\,|\right)^2\le C(a) E(a,\bfu).
\end{align}
If we take $\bfu=\bfG\bfdelta_{j_0}$ in \eqref{best-sob2}, 
then we have
\begin{align*}&
\left(\max_{0\le j\le 59}|\,{}^t\bfdelta_{j}\bfG\bfdelta_{j_0}\,|\right)^2 
\le C(a) E(a,\bfG\bfdelta_{j_0}) = C(a)^2.
\end{align*}
Combining this and a trivial inequality
\begin{align*}
C(a)^2=\left({}^t\bfdelta_{j_0}\bfG\bfdelta_{j_0}\right)^2
\le \left(\max_{0\le j\le 59}|\,{}^t\bfdelta_{j}\bfG\bfdelta_{j_0}\,|\right)^2,
\end{align*} 
we have
\begin{align*}
\left(\max_{0\le j\le 59}|\,{}^t\bfdelta_{j}\bfG\bfdelta_{j_0}\,|\right)^2 =
C(a) E(a,\bfG\bfdelta_{j_0})
.
\end{align*}
This shows that if we replace $C$ by $C(a)$ in an inequality \eqref{dsob2}, 
the equality holds for $\bfu=\bfG\bfdelta_{j_0}$. 
\eqref{C0a} follows from \eqref{C0alam} by simple calculation. 
\eqref{C0a} also follows from the observation of 
the diagonal value of $\bfG$ 
after calculating $\bfG$ by Mathematica. This completes the proof of Theorem \ref{th2}.
\hfill{$\blacksquare$}

\section*{Acknowledgment}
This research is supported by J. S. P. S. Grant-in-Aid for Scientific Research (C)
No.20540199, No.25400146, No.21540210 and No.26400175.


\end{document}